\documentclass{article}
\usepackage{amsfonts, amssymb, amsmath}
\usepackage[russian]{babel}
\begin{document}
\large
\begin{center}

{\bf\Large ОN A GENERALIZATION

OF BERNOULLI AND EULER NUMBERS}

{\bf A. Sarantsev}

\bigskip

{\it 10th International Seminar "Discrete Mathematics and its Applications"\ , 2010}

{\it University of Washington, PhD Student}

{\it E-mail: ansa1989@u.washington.edu}

\end{center}

\medskip

We shall consider permutations on the set $\{0, \ldots, N\}$ (i.e. one-to-one mappings from this set onto itself), where $N \in \mathbb Z_+$.
We say that such permutation $\tau$ {\it increases} ({\it decreases})
on the set $\{p, \ldots, q\} \subseteq \{0, \ldots, N\}$ if  for any $j = p, \ldots, q - 1$ we have
$\tau(j) < \tau(j + 1)$ ($\tau(j) > \tau(j + 1)$).

Take $n \in \mathbb N$,~ $i_1, j_1, \ldots, i_n, j_n \in \mathbb Z_+$,~ $N := i_1 + j_1 + \ldots + i_n + j_n$. The number of permutations on $\{0, \ldots, N\}$ which increase on $\{s''_{k - 1}, \ldots, s_k'\}$ and decrease on $\{s'_k, \ldots, s''_k\}$ for each $k = 1, \ldots, n$ will be denoted as $\Omega(i_1, j_1, \ldots, i_n, j_n)$. From now on $s'_k := i_1 + j_1 + \ldots + i_k$,~ $s''_k := i_1 + j_1 + \ldots + i_k + j_k$.

The last argument $j_n$ is omitted if $j_n = 0$. It is convenient to let $\Omega$ with no arguments be equal to $1$.

{\bf Lemma 1.} {\it 
$\Omega(p) = 1$,~ $\Omega(p, q) = (p+q)!/(p!q!)$ for $p, q \in \mathbb Z_+$.}

This is a generalization of {\it Andre problem}: to find $b_n := \Omega(1, \ldots, 1)$ ($n - 1$ unit arguments). Define {\it Bernoulli numbers} $B_n$ and {\it Euler numbers} $E_n$ as the coefficients in these series:
$$
\frac{x}{e^x - 1} = \sum\limits_{n = 0}^{+\infty}\frac{B_n}{n!}x^n,~~ \frac1{\ch x} = \sum\limits_{n = 0}^{+\infty}\frac{E_n}{n!}x^n.
$$
Then (see [1], Chapter 3, \$1)
$$
b_{2n - 1} = \frac{(-1)^{n-1}}n2^{2n - 1}(2^{2n} - 1)B_{2n},~~ b_{2n} = (-1)^nE_{2n}.
$$
Hence the $\Omega$ numbers are a generalization of Bernoulli and Euler numbers and binomial coefficients (recall Lemma 1).

{\bf Theorem 1.} {\it For any $n, i_1, \ldots, i_n \in \mathbb N$ we have $\Omega(i_1, \ldots, i_n) = \Omega(i_1 - 1, \ldots, i_n) + \Omega(i_1, i_2 - 1, i_3, \ldots, i_n) + \ldots + \Omega(i_1, \ldots, i_n - 1).$}



Consider a generating function
$$
F_n(x_1, \ldots, x_n) := \sum\limits_{i_1 = 0}^{+\infty}\ldots\sum\limits_{i_n = 0}^{+\infty}\Omega(i_1, \ldots, i_n)x_1^{i_1}\ldots x_n^{i_n}.
$$
It is not difficult to deduce $\Omega(i_1, \ldots, i_n) \le (i_1 + \ldots + i_n)!/(i_1!\ldots i_n!)$ from Theorem 1. Therefore, this series absolutely converges if $|x_1| + \ldots + |x_n| < 1$. For any $n$ the function $F_n$ is rational. One can calculate it for any given $n$. For example,
$$
F_3(x_1, x_2, x_3) = \frac{1 - x_1 - x_3}{(1 - x_1)(1 - x_3)(1 - x_1 - x_2 - x_3)},
$$
$$
F_4(x_1, x_2, x_3, x_4) = \frac1{1 - x_1 - x_2 - x_3 - x_4}\times
$$
$$
\times\biggl(\frac{(1 - x_1)(1 - x_1 - x_2 - x_4)}{(1 - x_1 - x_2)(1 - x_1 - x_4)} +
\frac{(1 - x_4)(1 - x_1 - x_3 - x_4)}{(1 - x_1 - x_4)(1 - x_3 - x_4)} - 1\biggr).
$$
The formulae for $\Omega$ with $3, 4$ or more arguments can be obtained form these expressions.

Let us introduce {\it exponential generating functions}: let $\Omega'(i_1, \ldots, i_n) := \Omega(i_1, \ldots, i_n)/(i_1 + \ldots + i_n + 1)!$ and
$$
G_n(x_1, \ldots, x_n) := \sum\limits_{i_1 = 0}^{+\infty}\ldots\sum\limits_{i_n = 0}^{+\infty}\Omega'(i_1, \ldots, i_n)x_1^{i_1}\ldots x_n^{i_n}.
$$
This term is commonly used for a generating function of $\Omega(i_1, \ldots, i_n)/(i_1!\ldots i_n!)$. But this is not convenient for us. - но это нам будет неудобно. Note that $\Omega'(i_1, \ldots, i_n)$ is a probability that an arbitrarily chosen permutation on $\{0, \ldots, i_1 + \ldots + i_n\}$ has the required monotonicity intervals.
The series for $G_n$ converges for any values of its arguments.

{\bf Theorem 2.} {\it For any $n, i_1, j_1, \ldots, i_n, j_n \in \mathbb N$, $N := i_1 + \ldots + j_n$ we have
$$
\Omega'(i_1, \ldots, j_n) = \frac1{N+1}\sum_{k=1}^n\Omega'(i_1, j_1, \ldots, i_k - 1)\Omega'(j_k - 1, \ldots, i_n, j_n),
$$
$$
\Omega'(i_1, \ldots, j_n) = \frac1{N + 1}\sum_{k=0}^n\Omega'(i_1, j_1, \ldots, j_k - 1)\Omega'(i_{k+1} - 1, \ldots, i_n, j_n).
$$}

It is straightforward to state and prove the similar theorems for odd number of arguments of $\Omega'$. The following theorem can be easily deduced from these four ones:

{\bf Theorem 3.} {\it For any $n, i_1, \ldots, i_n \in \mathbb N$, $N := i_1 + \ldots + i_n$ we have
$$
\Omega'(i_1, \ldots, i_n) = \frac1{2(N + 1)}\sum\limits_{k=0}^n\Omega'(i_1, \ldots, i_k - 1)\Omega'(i_{k+1} + 1, \ldots, i_n).
$$}
Let us apply this theory to probabilistic problems. Suppose $X_t,~ t \in \mathbb Z$ are independent identically distributed random variables with continuous distribution function.

{\bf Theorem 4.} {\it Let $n \in \mathbb N$, $i_1, j_1, \ldots, i_n, j_n \in \mathbb Z_+$, $N := i_1 + \ldots + j_n$. Тhen the following event: for all $k = 1, \ldots, n$
$$
X_{s''_{k-1}} < X_{s''_{k-1} + 1} < \ldots < X_{s'_k} > X_{s'_k + 1} > \ldots > X_{s''_k},
$$
has the probability $\Omega'(i_1, j_1, \ldots, i_n, j_n)$.}

Let us call any $n \in \mathbb Z$ {\it the maximum point} ({\it the minimum point}) if $X_n > X_{n - 1}$, $X_{n + 1}$ ($X_n < X_{n - 1}, X_{n + 1}$). The maximum points alternate with the minimum points. From now on we shall suppose that $0$ is a maximum point. Let $\mu_0$ be the distance from $0$ to the next minimum point, let $\mu_1$ be the distance from this point to the next maximum point, etc. Then the sequence $(\mu_t)$ is strictly stationary. It is easy to prove that for any $k_0, \ldots, k_n \in \mathbb N$~ $\mathbf P\{\mu_0 = k_0, \ldots, \mu_n = k_n\} = 3\Omega'(1, k_0, \ldots, k_n, 1)$.

After some calculations, we obtain $\mathbf P\{\mu_0 = k\} =  3(k^2 + 3k + 1)/(k + 3)!$,~ $\mathbf E\mu_0 = 3/2$, ~$Var\mu_0 = 6e - 63/4 \approx 0.560$, and $corr(\mu_0, \mu_1) = (2e^2 - 8e + 7)/(8e - 21) \approx 0. 0427$.

{\bf Conjecture.} The sequence $(\mu_t)$ satisfies the Law of Large Numbers:
$$
\frac1{n}\sum_{k=0}^{n-1}\mu_k \to \mathbf E\mu_0 = \frac32, ~ n \to +\infty.
$$
It is sufficient to prove that $(\mu_t)$ is ergodical. If this conjecture is true, one may try to construct a statistical test which checks whether any given observations $X_1, \ldots, X_{n}$ are independent and identically distributed.

\smallskip

\centerline{\bf References}

1. V. N. Sachkov, Introduction to Combinatorical Methods in Discrete Mathematics. Мoscow, MCCME Publishing House, 2004.

\end{document}